\documentclass[11pt]{article}
\usepackage{amsfonts}
\usepackage{amsmath}
\usepackage{amssymb}
\usepackage{hyperref}
\usepackage{tikz}

\newenvironment{acknowledgement}{\par\medskip\noindent\textbf{Acknowledgment.}\quad}{\par\medskip}

\def\email#1{\href{mailto:#1}{\texttt{#1}}}

\begin{document}

\title{An Improved Lower Bound for Maximal Partial Room Squares}
\author{
  Peter Horak\\[0.5em]
  University of Washington Tacoma\\
  \email{horak@uw.edu}
  \and
  Dongryul Kim\\[0.5em]
  Imperial College London\\
  \email{d.kim2@imperial.ac.uk}
}
\date{
  \textit{Dedicated to Alex Rosa on the occasion of his 90th birthday.}
}

\maketitle

\begin{abstract}
We establish a new asymptotic lower bound for the volume of a maximal
partial Room square of even order $n$. We prove that its volume is at least $%
\frac{2\sqrt{66}-6}{57}n^2-O(n)$, improving the leading coefficient $1/8$ in
the lower bound obtained by Meszka and Rosa. The proof combines two counting
arguments, one based on a transversal of empty cells and the other on the
graph of occurring symbol pairs.
\end{abstract}


\section{Introduction}

A \emph{Room square} of order $n$ is an $(n-1)\times(n-1)$ array on an $n$%
-element set $S$ in which each cell is either empty or contains an unordered
pair of distinct elements of $S$, subject to the following conditions: every
element of $S$ occurs exactly once in each row and in each column, and every
unordered pair of distinct elements of $S$ occurs in exactly one cell.
Necessarily, $n$ is even. The classical existence theorem of Mullin and
Wallis~\cite{MullinWallis1975} establishes the existence of Room squares of
every even order $n\geq 8$, while no Room squares of orders $4$ or $6$ exist.

A \emph{partial Room square} of order $n$ is defined by relaxing the last
two conditions: each element of $S$ occurs \emph{at most once} in each row
and each column, and each unordered pair occurs \emph{at most once} in the
array. Its \emph{volume}, denoted by $V$, is the number of occupied cells.
Thus a Room square has volume $\binom{n}{2}$. A partial Room square is \emph{%
maximal} if no unordered pair can be inserted into any empty cell without
violating one of these requirements. We refer to a maximal partial Room
square of order $n$ as an \emph{MPRS}$(n)$. Maximality is an inclusion
property and does not imply that the volume is as large as possible.
Throughout this paper, $n$ is even.

The problem of determining the minimum number of occupied cells in a maximal
partial combinatorial array has a longer history. For partial Latin squares,
Horak and Rosa~\cite{HorakRosa1993} proved that a maximal partial Latin
square of order $n$ has at least $\lceil n^{2}/2\rceil $ occupied cells, and
that this bound is attained for every $n$. They also determined the spectrum
of possible numbers of occupied cells up to a linear number of values.

A natural extension concerns pairs of orthogonal partial Latin squares. In
this setting, maximality means that no ordered pair of symbols can be
inserted into an empty cell while preserving the Latin property in both
squares and their orthogonality. Donovan, Grannell and Yazici~\cite%
{DonovanGrannellYazici2026} proved that every maximal pair of orthogonal partial
Latin squares of order $n$ has at least $n^2/3$ occupied cells. They also
showed that, for $n\geq 21$, the minimum is exactly $\lceil n^2/3\rceil$ and
characterized the extremal arrays. The same lower bound was obtained
independently by Horak, Meszka and Rosa in a manuscript currently being
prepared for submission; their work also investigates the spectrum of
volumes of maximal orthogonal partial Latin squares~\cite%
{HorakMeszkaRosaForthcoming}.

These results provide a natural context for the corresponding question for
partial Room squares, where each occupied cell contains an unordered pair of
symbols rather than a single symbol or an ordered pair.

An important structural difference between these problems is the role of
conjugacy. A partial Latin square can be represented as a set of ordered
triples, and a pair of orthogonal partial Latin squares as a set of ordered
quadruples. In both settings, permutations of the coordinates yield
conjugate arrays of the same kind. Thus rows, columns, and symbol
coordinates play symmetric roles: arguments involving the minimum frequency
of one coordinate can be transferred to any other coordinate by conjugation.
For a partial Room square, however, the symbols occur as unordered pairs in
occupied cells, and the symbol coordinates are not interchangeable with the
row and column coordinates. In particular, the minimum number of occurrences
of a symbol and the minimum number of occupied cells in a row or column are
distinct parameters, with no conjugacy relation forcing them to agree. This
asymmetry makes the analysis of maximal partial Room squares more difficult.
Our proof reflects this distinction by using separate parameters for the two
types of minimum frequency.

The possible volumes of maximal partial Room squares were studied by Meszka
and Rosa~\cite{MeszkaRosa2021}. Writing 
\begin{equation*}
M(n)=\{V:\text{there exists an MPRS}(n)\text{ of volume }V\},\quad
m(n)=\min M(n),
\end{equation*}%
they investigated the \emph{volume spectrum} $M(n)$, with particular
emphasis on its lower endpoint. Their principal result shows that, for $%
n=2m\geq 14$, 
\begin{equation}
\begin{cases}
\lbrack m^{2}-m,\,\binom{n}{2}]\cap \mathbb{Z}\subseteq M(n), & \text{if $m$
is even}, \\[2pt] 
\lbrack m^{2}-m+1,\,\binom{n}{2}]\cap \mathbb{Z}\subseteq M(n), & \text{if $%
m $ is odd}.%
\end{cases}
\label{A}
\end{equation}%
They conjectured that the respective left endpoints in \eqref{A} are the
actual values of $m(n)$. Thus, if this conjecture holds, their theorem
determines the entire spectrum $M(n)$. Their results establish the existence
of maximal partial Room squares for every volume in the indicated intervals,
thereby providing an extensive description of the volume spectrum.

For perspective, an MPRS of volume exactly $n^{2}/4=m^{2}$ has a simple
construction for every $n=2m$ with even $m\geq 8$. Partition $S$ into
disjoint $m$-element sets $A=\{a_{1},\ldots ,a_{m}\}$ and $B=\{b_{1},\ldots
,b_{m}\}$. In the first $m-1$ columns of a $(2m-1)\times (2m-1)$ array,
place a Room square of order $m$ on $A$ in the first $m-1$ rows and another
Room square of order $m$ on $B$ in the next $m-1$ rows. The cells in the
last row of the first $m-1$ columns are empty. Further, all cells in the
last $m$ columns are empty except for those in the last row, which are
filled with $m$ pairs $\{a_{i},b_{i}\}$, $1\leq i\leq m$. The resulting
partial Room square is maximal, and its volume is 
\begin{equation*}
V=2\binom{m}{2}+m=m^{2}=\frac{n^{2}}{4}.
\end{equation*}%

\begin{figure}[ht]
  \centering
  \begin{tikzpicture}
    \draw[step=0.5] (1.9,0) grid (4.5,0.5);
    \draw[thick] (0,0) rectangle (4.5,4.5);
    \draw[thick] (0,0.5) rectangle (2,2.5);
    \draw[thick] (0,2.5) rectangle (2,4.5);
    \node at (1,3.5) {$RS_A$};
    \node at (1,1.5) {$RS_B$};
    \draw[->] (2.5,1) node[above] {$a_1b_1$} -- (2.25,0.25);
    \draw[->] (4,1) node[above] {$a_mb_m$} -- (4.25,0.25);
    \node at (3.285,0.25) {$\cdots$};
  \end{tikzpicture}
  \caption{A drawing of the construction of a MPRS of volume $n^2/4$}
\end{figure}
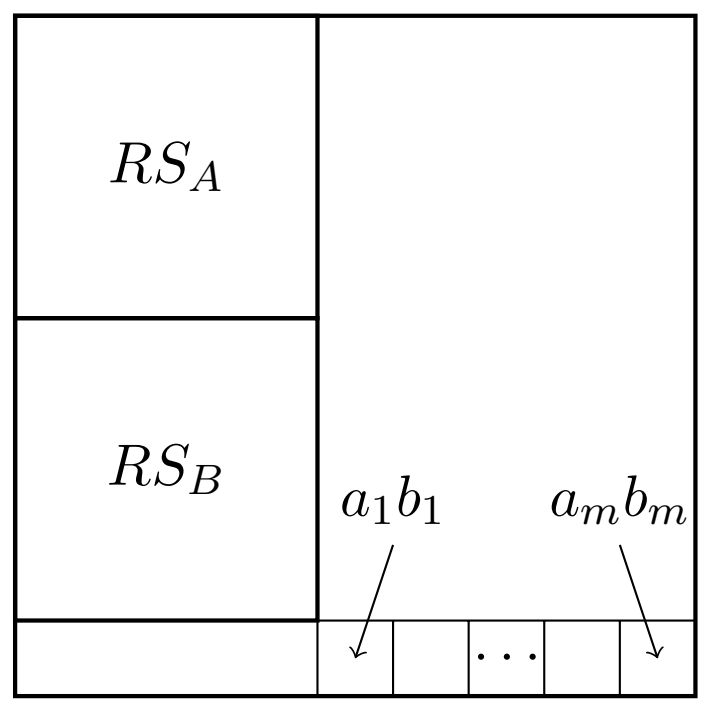

On the other hand, the general lower bound proved by Meszka and Rosa~\cite[%
Lemma~29]{MeszkaRosa2021} states that every MPRS of even order $n\geq 6$
satisfies 
\begin{equation*}
V>\frac{n(n-1)}{8}.
\end{equation*}%
The lower bound of Meszka and Rosa and their conjectured minimum leave an
interesting interval of possible values for $m(n)$. Motivated by their work,
we investigate whether the lower bound can be strengthened. We establish the
asymptotic lower bound 
\begin{equation*}
m(n)\geq \frac{2\sqrt{66}-6}{57}\,n^{2}-O(n)\qquad (n\rightarrow \infty ,\ n%
\text{ even}).
\end{equation*}%
Thus, the leading coefficient in the known lower bound increases from $1/8$
to $0.1797908212\ldots $. Our result constitutes a further step toward
determining the minimum volume conjectured by Meszka and Rosa.

\section{Two elementary bounds}

Let $RS$ be an MPRS$(n)$ with volume $V$. Let $m_{1}$ denote the minimum
frequency of a symbol, and let $m_{2}$ denote the minimum number of occupied
cells among all rows and columns.

Since each occupied cell contains two symbols and there are $n-1$ rows, we
immediately obtain 
\begin{equation}
V\geq\frac{nm_{1}}{2},\qquad V\geq(n-1)m_{2}.  \label{1}
\end{equation}

\section{A bound based on a transversal of empty cells}

Choose a symbol $x$ occurring exactly $m_{1}$ times, in the pairs $%
xy_{1},\ldots,xy_{m_{1}}$. These pairs have distinct rows and distinct
columns. By permuting the rows and columns, we may assume that they occupy
the main diagonal of the upper-left $m_{1}\times m_{1}$ block $A$. Partition
the array into nine blocks as follows: 
\begin{equation*}
RS=\left( 
\begin{array}{ccc}
A & B & G \\ 
C & D & H \\ 
I & J & K%
\end{array}
\right) .
\end{equation*}
For any block $X$, write $|X|$ for its number of occupied cells.

Let $p_{1}$ be the \emph{maximum cardinality} of a transversal of empty
cells in the subarray formed by the rows and columns outside $A$. Permuting
these remaining rows and columns, we may place the transversal of empty
cells $P_{1}$ on the main diagonal of $D$, where $D$ is $p_{1}\times p_{1}$.
Write $q_{1}$ for the side length of $K$. Thus 
\begin{equation}
m_{1}+p_{1}+q_{1}=n-1.  \label{2}
\end{equation}%
The first block row has $m_{1}$ rows, the second has $p_{1}$ rows, and the
third has $q_{1}$ rows; the corresponding block columns have the same sizes.

Every cell of $K$ is occupied, for otherwise an empty cell of $K$ could be
added to $P_{1}$. Consequently, 
\begin{equation}
|K|=q_{1}^{2}.  \label{3}
\end{equation}%
Moreover, for each $i\in \{1,\ldots ,p_{1}\}$, either the $i$th row of $H$
is completely occupied or the $i$th column of $J$ is completely occupied.
Indeed, if both contained an empty cell, we could delete the empty diagonal
cell $(i,i)$ of $D$ from $P_{1}$ and replace it by those two empty cells,
obtaining a transversal of empty cells of size $p_{1}+1$. Hence 
\begin{equation}
|H|+|J|\geq p_{1}q_{1}.  \label{4}
\end{equation}

Every row and every column has at least $m_{2}$ occupied cells. Summing over
the $m_{1}$ rows and the $m_{1}$ columns meeting $A$ gives 
\begin{equation}
|A|+|B|+|G|\geq m_{1}m_{2},\qquad|A|+|C|+|I|\geq m_{1}m_{2}.  \label{5}
\end{equation}

Now set 
\begin{equation*}
T_{1}=S\setminus\{x,y_{1},\ldots,y_{m_{1}}\}, \qquad|T_{1}|=n-1-m_{1}.
\end{equation*}
For every $t\in T_{1}$, the pair $xt$ does not occur in $RS$. In addition, $%
x $ occurs in none of the rows or columns meeting $D$, because all
occurrences of $x$ lie in $A$. Thus, for any empty cell $(i,j)$ of $D$,
maximality forces $t$ to occur in row $i$ or column $j$; otherwise we could
insert $xt$ into that cell.

Apply this observation to the $p_{1}$ empty diagonal cells of $D$. Each
symbol $t\in T_{1}$ must occur in the row or column of each such cell, so
altogether there are at least $(n-1-m_{1})p_{1}$ symbol--diagonal incidences
to account for. An occupied cell of $B$, $C$, $H$, or $J$ accounts for at
most two incidences, while an occupied cell of $D$ accounts for at most
four. Cells in the other blocks account for none. Therefore 
\begin{equation}
|B|+|C|+|H|+|J|+2|D|\geq\frac{(n-1-m_{1})p_{1}}{2}.  \label{6}
\end{equation}

Using the identity 
\begin{equation*}
\begin{array}{rcl}
2V & = & (|A|+|B|+|G|)+(|A|+|C|+|I|) \\ 
&  & +(|B|+|C|+|H|+|J|+2|D|) \\ 
&  & +(|H|+|J|)+2|K|+(|G|+|I|),%
\end{array}%
\end{equation*}
and combining \eqref{3}--\eqref{6} we obtain 
\begin{equation}
V\geq m_{1}m_{2}+\frac{(n-1-m_{1})p_{1}}{4}+\frac{p_{1}q_{1}}{2}+q_{1}^{2}.
\label{7}
\end{equation}

\section{A bound based on the graph of occurring pairs}

Define a simple graph $\Gamma$ with vertex set $S$, in which two distinct
symbols are adjacent if and only if their pair occurs in $RS$. Each occupied
cell corresponds to one edge, and the degree of a vertex is the frequency of
the corresponding symbol. In particular, 
\begin{equation}
2V=\sum_{v\in S}\deg_{\Gamma}(v).  \label{8}
\end{equation}

Choose a row $R$ containing exactly $m_{2}$ occupied cells. Its pairs form a
matching $M$ in $\Gamma$. Let $U$ be the set of the $2m_{2}$ vertices
covered by $M$, and put $T_{2}=S\setminus U$. Let 
\begin{equation*}
P_{2}=\{\{a_{1},b_{1}\},\ldots,\{a_{p_{2}},b_{p_{2}}\}\}
\end{equation*}
be a maximum matching of nonedges of $\Gamma$ with endpoints in $T_{2}$.
Write $W=\{a_{1},b_{1},\ldots,a_{p_{2}},b_{p_{2}}\}$ and $Q=T_{2}\setminus W$%
, with $|Q|=q_{2}$. Then 
\begin{equation}
2m_{2}+2p_{2}+q_{2}=n.  \label{9}
\end{equation}
Since $n$ is even, $q_{2}$ is even as well. Thus either $q_{2}=0$ or $%
q_{2}\geq2$. Since $P_{2}$ is a maximum matching of nonedges in $T_{2}$,
every two distinct vertices of $Q$ must be adjacent in $\Gamma$. Thus the
graph induced by $Q$ is complete.

For each pair $\{a_{i},b_{i}\}\in P_{2}$, neither symbol occurs in row $R$
and their pair is absent from $RS$. The row $R$ has $n-1-m_{2}$ empty cells.
By maximality, at least one of $a_{i},b_{i}$ must occur in the column of
each of these cells; otherwise their pair could be inserted. It follows that 
\begin{equation}
\sum_{w\in W}\deg_{\Gamma}(w)\geq p_{2}(n-1-m_{2}).  \label{10}
\end{equation}
Also, by the definition of $m_{1}$, 
\begin{equation}
\sum_{u\in U}\deg_{\Gamma}(u)\geq2m_{2}m_{1}.  \label{11}
\end{equation}

Suppose that $q_{2}\geq2$. Fix a pair $\{a,b\}\in P_{2}$. If each vertex of $%
Q$ is adjacent to at least one of $a,b$, there are at least $q_{2}$ edges
between $\{a,b\}$ and $Q$. Otherwise, choose $z\in Q$ that is adjacent to
neither $a$ nor $b$. Every vertex $z^{\prime}\in Q\setminus\{z\}$ must then
be adjacent to both $a$ and $b$: if, for example, $az^{\prime}$ were also a
nonedge, we could replace the nonedge $ab$ in $P_{2}$ by the two disjoint
nonedges $az^{\prime}$ and $bz$. The same argument applies with $a$ and $b$
interchanged. In this case there are at least $2(q_{2}-1)\geq q_{2}$ edges
from $\{a,b\}$ to $Q$. Summing over the pairs of $P_{2}$ shows that at least 
$p_{2}q_{2}$ edges join $W$ to $Q$.

Since $Q$ induces a complete graph, for $q_{2}\geq2$ we have 
\begin{equation}
\sum_{v\in Q}\deg_{\Gamma}(v)\geq q_{2}(q_{2}-1)+p_{2}q_{2}.  \label{12}
\end{equation}
The same inequality is trivially true when $q_{2}=0$. These are the only
possible cases. As $S$ is the disjoint union of $U$, $W$, and $Q$, equations
\eqref{8}, \eqref{10}, \eqref{11}, and \eqref{12} give 
\begin{equation}
V\geq m_{2}m_{1}+\frac{p_{2}(n-1-m_{2})}{2}+\frac{q_{2}(q_{2}-1+p_{2})}{2}.
\label{13}
\end{equation}

\section{Minimization of the lower bounds}

We now justify the global minimization of the necessary inequalities, rather
than merely verifying a candidate choice of parameters. Put 
\begin{equation*}
c=\frac{2\sqrt{66}-6}{57}=\frac{2}{3+\sqrt{66}}=0.1797908212\ldots .
\end{equation*}

\medskip\noindent\textbf{Theorem 1.} \textit{Every maximal partial Room
square on an even set of $n$ symbols has volume} 
\begin{equation*}
V\geq cn^{2}-O(n).
\end{equation*}%
\noindent\textit{Moreover, $c$ is the best leading coefficient obtainable by
minimizing simultaneously the necessary bounds \eqref{1}, \eqref{7}, and
\eqref{13}, subject to \eqref{2} and \eqref{9}.}

\medskip \noindent \textit{Proof.} Write $a=m_{1}/n$ and $b=m_{2}/n$. Thus $%
0\leq a\leq 1$ and $0\leq b\leq 1/2$. We first eliminate the auxiliary
variables. Setting $t=(n-1-m_{1})/n$ and $u=q_{1}/n$, the relation $%
p_{1}=n-1-m_{1}-q_{1}$ gives the exact identity 
\begin{equation*}
\frac{(n-1-m_{1})p_{1}}{4n^{2}}+\frac{p_{1}q_{1}}{2n^{2}}+\frac{q_{1}^{2}}{%
n^{2}}=\frac{t^{2}}{4}+\frac{tu}{4}+\frac{u^{2}}{2}\geq \frac{t^{2}}{4}.
\end{equation*}%
Consequently, the transversal of empty cells estimate implies, uniformly
over all admissible parameters, 
\begin{equation}
\frac{V}{n^{2}}\geq ab+\frac{(1-a)^{2}}{4}-O(n^{-1}).  \label{eq:reducedT}
\end{equation}

For the graph estimate put $r=q_{2}/n$ and use $p_{2}/n=(1-2b-r)/2$, where $%
0\leq r\leq 1-2b$. Direct substitution into \eqref{13} yields 
\begin{equation*}
\frac{V}{n^{2}}\geq ab+\frac{(1-b)(1-2b)}{4}-\frac{br}{4}+\frac{r^{2}}{4}%
-O(n^{-1}).
\end{equation*}%
The quadratic in $r$ has its unconstrained minimum at $r=b/2$. In
particular, if $b\leq 2/5$, this value is feasible and hence 
\begin{equation}
\frac{V}{n^{2}}\geq ab+g(b)-O(n^{-1}),\qquad g(b)=\frac{1-3b+\frac{7}{4}b^{2}%
}{4}.  \label{eq:reducedG}
\end{equation}%
We also have from \eqref{1} 
\begin{equation}
\frac{V}{n^{2}}\geq b-O(n^{-1}).  \label{eq:reducedB}
\end{equation}%
If $b\geq c$, \eqref{eq:reducedB} proves the required bound. It remains to
consider $0\leq b\leq c<2/5$, where both \eqref{eq:reducedT} and %
\eqref{eq:reducedG} apply.

Let $D(b)=1-3b+\frac74 b^2$ and $A(b)=1-\sqrt{D(b)}$. For $0\leq b\leq c$,
one has $0\leq A(b)<1-2b$; the latter inequality follows from $D(b)>4b^2$ on
this interval. The two functions 
\begin{equation*}
T(a,b)=ab+\frac{(1-a)^2}{4},\qquad G(a,b)=ab+g(b)
\end{equation*}
meet at $a=A(b)$ (the other algebraic solution lies outside $[0,1]$). For
fixed $b$, $T$ decreases on $[0,A(b)]$, while $G$ increases with $a$ on $%
[A(b),1]$. It follows that 
\begin{equation}  \label{eq:fixedb}
\min_{0\leq a\leq1}\max\{T(a,b),G(a,b)\} = F(b):=b\bigl(1-\sqrt{D(b)}\bigr)+%
\frac{D(b)}4.
\end{equation}
An elementary differentiation shows that $F$ is decreasing on $[0,c]$.
Indeed, 
\begin{equation*}
F^{\prime }(b)=1-\sqrt{D(b)}+\frac{7b-6}{8} +\frac{b(6-7b)}{4\sqrt{D(b)}}<0
\qquad (0\leq b\leq c).
\end{equation*}
The final inequality follows directly upon substituting $D(b)=1-3b+7b^2/4$
and using $0\leq b\leq c<9/50$ (or by squaring the two positive sides after
rearrangement). Thus $F(b)\geq F(c)$. The equation $57c^2+12c-4=0$ gives $%
D(c)=16c^2$, and hence 
\begin{equation*}
F(c)=c(1-4c)+4c^2=c.
\end{equation*}
Combining \eqref{eq:reducedT}--\eqref{eq:fixedb} proves $V/n^2\geq
c-O(n^{-1})$.

To see that these inequalities alone yield no larger leading coefficient,
take admissible integers with 
\begin{align*}
m_1&=(1-4c)n+O(1), & p_1&=4cn+O(1), & q_1&=0, \\
m_2&=cn+O(1), & p_2&=(\tfrac12-\tfrac{5c}{4})n+O(1), & q_2&=(\tfrac
c2)n+O(1).
\end{align*}
Choose $m_1,m_2$ by rounding, set $p_1=n-1-m_1$, choose $q_2$ as a nearby
even nonnegative integer, and set $p_2=(n-2m_2-q_2)/2$. All these quantities
are nonnegative for sufficiently large even $n$, and the partition
identities hold. Substitution shows that the leading coefficients of the
four right-hand sides in \eqref{1}, \eqref{7}, and \eqref{13}, are,
respectively, 
\begin{equation*}
\frac{1-4c}{2}<c,\qquad c,\qquad c(1-4c)+4c^2=c,\qquad c(1-4c)+\frac{%
1-3c+\frac74c^2}{4}=c.
\end{equation*}
This establishes optimality for the system of necessary inequalities. It
does not assert the existence of maximal partial Room squares attaining the
bound. \hfill $\square$

\medskip

\section{Improvement of the lower bound}

Meszka and Rosa~\cite{MeszkaRosa2021} established that every maximal partial
Room square of even order $n\geq 6$ has volume greater than $n(n-1)/8$. The
inequalities obtained in the present paper, together with their
optimization, yield the stronger asymptotic lower bound 
\begin{equation*}
V\geq \frac{2\sqrt{66}-6}{57}\,n^{2}-O(n).
\end{equation*}%
Consequently, the leading coefficient in the lower bound increases from $1/8$
to $0.1797908212\ldots $. Although the conjectured minimum volume remains
undetermined, this result narrows the gap between the known lower bound and
the conjecture of Meszka and Rosa.

\begin{acknowledgement}
The authors acknowledge that the calculations in Section 5 and the
construction of a maximal partial Room square of volume $\frac{n^{2}}{4}$
presented in the Introduction were produced by ChatGPT Plus (OpenAI). All
other mathematical results in this paper were obtained independently by the
authors.
\end{acknowledgement}

\end{document}